\documentclass[11pt]{article}
\usepackage{amsmath, amscd, amssymb, mathrsfs, latexsym, epsfig, color, amsthm, mathtools, tikz-cd, array, adjustbox, bbm, tikz, enumitem, relsize, xcolor}

\usepackage{titling}

\tikzset{set/.style={draw,ellipse,inner sep=0pt,align=left, anchor=west}}


\DeclarePairedDelimiter\ceil{\lceil}{\rceil}
\DeclarePairedDelimiter\floor{\lfloor}{\rfloor}
\usepackage[normalem]{ulem}

\makeatletter
\def\moverlay{\mathpalette\mov@rlay}
\def\mov@rlay#1#2{\leavevmode\vtop{%
		\baselineskip\z@skip \lineskiplimit-\maxdimen
		\ialign{\hfil$\m@th#1##$\hfil\cr#2\crcr}}}
\newcommand{\charfusion}[3][\mathord]{
	#1{\ifx#1\mathop\vphantom{#2}\fi
		\mathpalette\mov@rlay{#2\cr#3}
	}
	\ifx#1\mathop\expandafter\displaylimits\fi}
\makeatother

\numberwithin{equation}{section}
\newtheorem{theorem}{Theorem}[section]
\newtheorem{proposition}[theorem]{Proposition}

\newtheorem{corollary}[theorem]{Corollary}

\newtheorem{lemma}[theorem]{Lemma}

\theoremstyle{definition}

\DeclareMathOperator\lk{\mathrm{lk}}

\DeclareMathOperator\Stack{\mathrm{Stack}}

\newcommand{\Pyr}{\operatorname{Pyr}}

\tikzset{
	labl/.style={anchor=south, rotate=90, inner sep=.5mm}
}


\usepackage[normalem]{ulem}
\usepackage{ifxetex,ifluatex,stmaryrd}
\newif\ifxetexorluatex
\ifxetex
  \xetexorluatextrue
\else
  \ifluatex
    \xetexorluatextrue
  \else
    \xetexorluatexfalse
  \fi
\fi

\ifxetexorluatex
  \usepackage{fontspec}
\else
  \usepackage[T1]{fontenc}
  \usepackage[utf8]{inputenc}
  \usepackage{lmodern}
\fi

\usepackage{hyperref}
\usepackage{amsmath,amsthm,amsfonts,amssymb,graphicx,mathtools,array,bbm,soul,xcolor}
\usepackage{subcaption}
\usepackage[margin=1in]{geometry}

\usepackage[normalem]{ulem}

\theoremstyle{definition}

\title{A Lower Bound Theorem for strongly regular CW spheres with up to $2d+1$ vertices}
\author{Lei Xue}

\begin{document}
	\maketitle
	\begin{abstract}
		In 1967, Gr\"unmbaum conjectured that any $d$-dimensional polytope with $d+s\leq 2d$ vertices has at least
		\[\phi_k(d+s,d) = {d+1 \choose k+1 }+{d \choose k+1 }-{d+1-s \choose k+1 } \]
		$k$-faces. This conjecture along with the characterization of equality cases was recently proved by the author \cite{myLBTpaper}. In this paper, several extensions of this result are established. Specifically, it is proved that lattices with the diamond property (for example, abstract polytopes) and $d+s\leq 2d$ atoms have at least $\phi_k(d+s,d)$ elements of rank $k+1$. Furthermore, in the case of face lattices of strongly regular CW complexes representing normal pseudomanifolds with up to $2d$ vertices, a characterization of equality cases is given. Finally, sharp lower bounds on the number of $k$-faces of strongly regular CW complexes representing normal pseudomanifolds with $2d+1$ vertices are obtained. These bounds are given by the face numbers of certain polytopes with $2d+1$ vertices.
	\end{abstract}
	
	\section{Introduction}
	This paper is devoted to the study of the minimal face numbers of strongly regular CW complexes representing normal pseudomanifolds. A great deal of research has been done on the face numbers of polytopes in the past fifty years. In 1970, McMullen \cite{MR283691} established the Upper Bound Theorem, providing tight upper bounds on the number of $k$-faces a $d$-dimensional polytope with $n$ vertices can have. Later, Barnette (see \cite{MR298553}, \cite{MR360364}, and \cite{MR328773}) proved the Lower Bound Theorem for {\em simplicial} polytopes; his result provides tight lower bounds on the number of $k$-faces a $d$-dimensional simplicial polytope with $n$ vertices can have. Moreover, the lower bound theorem was proved to hold even in the generality of simplicial complexes that are normal pseudomanifolds (see \cite{MR877009}, \cite{MR1314963}, and \cite{MR2636638}). In 1980, Billera and Lee \cite{MR551759} and Stanley \cite{MR563925} completely characterized the face numbers of all simplicial (and by duality also simple) polytopes. Their result is known as the $g$-theorem.  Billera and Lee \cite{MR551759} established sufficiency of the conditions while Stanley \cite{MR563925} proved their necessity. 
	
	Despite many advances on face numbers of simplicial polytopes, much less is known about the face numbers of general polytopes. In fact, the only plausible conjecture providing the lower bounds on face numbers was made by Gr\"unbaum in \cite[p.~184]{Gr1-2} (see also \cite[p.~265]{MR250188}). He conjectured that for general $d$-dimensional polytopes with $d+s\leq 2d$ vertices, the number of $k$-faces is at least
	\[\phi_k(d+s,d) = {d+1 \choose k+1 }+{d \choose k+1 }-{d+1-s \choose k+1 }. \] He proved this conjecture for the cases of $s=2,\;3$, and $4$. The conjecture remained completely open for $s\geq 5$  until Pineda-Villavicencio, Ugon and Yost \cite{MR3899083} proved this conjecture for the number of edges, i.e., they verified the $k=1$ case. In 2020, we proved (see \cite{myLBTpaper}) the conjecture in full generality and characterized equality cases. 
	
	In this paper we extend our results to the following much more general classes of posets and complexes. For all undefined terminology and notation; see Section 2.
	\begin{itemize}
		\item We prove that Gr\"unbaum's conjecture holds in the generality of lattices with the diamond property. In other words, for any diamond lattice $L$ of rank $d+1$ with $d+s \leq 2d$ atoms, the number of elements of $L$ of rank $k+1$ is at least $\phi_{k}(d+s,d)$; see Theorem \ref{thm: inequality lattice}.
		\item For the case of face lattices of strongly regular CW complexes representing normal $(d-1)$-pseudomanifolds with up to $2d$ vertices, we characterize the equality cases; see Theorem \ref{thm:equality sphere}. Along the way, we show that any strongly regular normal $(d-1)$-pseudomanifold with $d+2$ facets is combinatorially equivalent to the boundary complex of a $d$-polytope; see Theorem \ref{thm: spheres are polytopes}.  
		
		\item Finally, in the class of strongly regular normal $(d-1)$-pseudomanifolds with $2d+1$ vertices, we determine the minimum value of the number of $k$-dimensional faces for all $1\leq k\leq d-1$. These minimum values are achieved by certain polytopes with $d+2$ and $d+3$ facets. See Section \ref{section: 2d+1 verts}.  The lower bound for these pseudomanifolds has two parts. Specifically, let $P$ be such a pseudomanifold. Then the following holds.
		\subitem If $P$ has $d+2$ facets, then $P$ is combinatorially equivalent to the boundary complex of a $d$-polytope (see Theorem \ref{thm: spheres are polytopes}), and for $1\leq k\leq d-2$,
		\begin{align}\label{eqn: 2d+1, d+2}
		f_k(P) \geq {d+1 \choose k+1} + {d \choose k+1} + {d-1 \choose k+1} - {\ceil{\frac{d}{2}} \choose k+1} - {\ceil{\frac{d}{2}} -1 \choose k+1}, \text{ see Lemma \ref{lem: bound for (2d+1, d+2)}.} 
		\end{align}
		
		\subitem If $P$ has at least $d+3$ facets, then for $1\leq k\leq d-2$,
		\begin{align}\label{eqn: 2d+1, d+3}
		f_k(P) \geq {d+1 \choose k+1} + {d \choose k+1} + {d-1 \choose k}, \text{ see Theorem \ref{thm: 2d+1}.}
		\end{align}
		
		Some of these results, but only for the case of {\em polytopes} with $2d+1$ vertices, were obtained by Pineda-Villavicencio and Yost in \cite{pvy2021}. At the end of \cite{pvy2021} they conjectured (based on computational data up to $d=100$) that for any $d$-polytope with $2d+1$ vertices, the lower bound for the numbers of $m$-dimensional faces has the following two parts:
		\begin{itemize}
			\item[(i)] When $1\leq m\leq \ceil{d/2}-1$, the bound is (\ref{eqn: 2d+1, d+3}) 
			\item[(ii)] When $m$ is larger, the lower bound comes from some other polytope with $d+2$ facets. 
		\end{itemize}
		Part (i) of this conjecture could be easily confirmed by comparing the two formulas (\ref{eqn: 2d+1, d+2}) and (\ref{eqn: 2d+1, d+3}) (see Appendix).

	\end{itemize}

	\section{Basics on polytopes and strongly regular CW complexes}
	
	We will first recall some definitions and introduce some notation. We refer the reader to Gr\"unbaum's and Ziegler's books (see \cite{Gr1-2} and \cite{MR1311028}) for all undefined notions related to polytopes, and to Bj\"orner's survey paper \cite{MR1373690} for notions related to CW complexes. By a {\bf polytope} we mean the convex hull of finitely many points in $\mathbb{R}^d$. A {\bf $d$-simplex}, denoted as $\Delta^d$, is the convex hull of $d+1$ affinely independent points. A {\bf face} of a polytope $P$ is the intersection of $P$ with a supporting hyperplane. It is known that a face of a polytope is a polytope. The dimension of a polytope is the dimension of its affine span. For brevity, we refer to a $k$-dimensional face as a $k$-face and to a $d$-dimensional polytope as a $d$-polytope. The $0$-faces are called vertices. The $1$-faces are called edges. The $(d-1)$-faces of a $d$-polytope are called {\bf facets}. We denote by $f_k(P)$ the number of $k$-faces of a polytope $P$.
	
	Let $P\subset \mathbb{R}^d$ be a $d$-polytope and $v$ a vertex of $P$. The {\bf vertex figure of $P$ at $v$}, $P/v$, is obtained by intersecting $P$ with a hyperplane $H$ that separates $v$ from the rest of the vertices of $P$. One property of vertex figures that will be very useful for us is that $(k-1)$-faces of $P/v$ are in bijection with $k$-faces of $P$ that contain $v$. If a vertex $v\in P$ is contained in exactly $d$ edges, then $v$ is a {\bf simple} vertex, otherwise $v$ is in more than $d$ edges and it is called {\bf nonsimple}.
	
	A {\bf CW complex} is any space $X$ which can be built in the following ``hierarchical'' way: start with a (finite) collection of isolated points called $X^0$, then attach one-dimensional disks $D^1$ to $X^0$ along the boundaries $S^0$, denote this object by $X^1$, then attach two-dimensional disks $D^2$ to $X^1$ along their boundaries $S^1$, writing $X^2$ for the new space, and so on, giving spaces $X^k$ for every $k$. A CW complex is any space that has this sort of decomposition into subspaces $X^k$. (The $X^k$s must exhaust all of $X$.) The (closed) $k$-disks used to build $X$ are called the $k$-dimensional {\bf faces}\footnote{The disks here are usually called (closed) {\bf cells} in the literature. We call them faces in this paper to be consistent with the notion of faces of polytopes.} or ($k$-faces). A CW complex $X$ is {\bf pure} if all of it facets have the same dimension.
	
	A CW complex $X$ is a {\bf CW sphere} if it is homeomorphic to a sphere. And it is {\bf regular} if every attaching map is a homeomorphism. Moreover, we call a CW sphere {\bf strongly regular} if, in addition, the intersection of any two faces is a face (possibly empty). For example, the boundary complex of a $d$-polytope is a strongly regular CW $(d-1)$-sphere. Also, any oriented matroid sphere (as defined in \cite{MR889977}) is a strongly regular CW sphere.
	
	Let $L = (L,\leq)$ be a poset (partially ordered set). Any (totally) ordered subset $t_0<t_1<\dots<t_{i}$ is a {\bf chain} in $L$ of length $i$. A finite poset is {\bf graded} if all the maximal chains have the same length. Throughout this paper, every poset we discuss will be finite and graded. Every poset will also be {\bf bounded}, i.e., have unique bottom and top elements $\hat{0}$ and $\hat{1}$. Let $\rho:L\to \mathbb{N}$ be the rank function. For the bottom element, the rank is $\rho (\hat{0}) = 0$. The rank of $L$, $\rho(L)$, is defined as $\rho(\hat{1})$. For $x,y$ two elements of $L$ where $x\leq y$, the {\bf interval} $[x,y]$ is the sub-poset of $L$ that contains all elements $z \in L$ such that $x\leq z\leq y$ (in particular, we can write $L$ itself as the interval $[\hat{0}, \hat{1}]$). When $x < y$, we say $x$ is {\bf below} $y$ and $y$ is {\bf above} or {\bf contain} $x$. If in addition there exists no $z \in L$ such that $x < z < y$, then $y$ {\bf covers} $x$. A poset is a {\bf lattice} if for every pair of elements, both the least upper bound and the greatest lower bound exist.

	If a graded lattice $L= [\hat{0}, \hat{1}]$ satisfies the condition that all rank-$2$ intervals are of the diamond shape (i.e., have four elements), then we call $L$ a {\bf diamond lattice} or {\bf $\Diamond$-lattice}. In a lattice of rank $d$, the rank-$1$ elements are called the {\bf atoms}, while the rank-$(d-1)$ elements are called the {\bf coatoms}. A lattice $L$ is called {\bf coatom-distinguishable} if given any two elements $t, s\in L\backslash \{\hat{0}, \hat{1}\}$ of the same rank, there exists a coatom $F\in L$ such that $t \leq F$ while $s\nleq F$. We will show that $\Diamond$-lattices are coatom-distinguishable (see Proposition \ref{prop: diamond implies coatom-dist}). 
	
	The {\bf open face poset} $L^o(X)$ of a CW complex $X$ is the collection of its proper faces ordered by inclusion, together with the unique bottom element $\hat{0}$. The {\bf face poset} $L(X)$ is $L^o(X)$ to which the unique top element $\hat{1}$ was added. If $X$ is a regular CW complex, then (1) the face poset $L(X)$ is graded, and (2) for every $x\in X$, the interval $[\hat{0},x]$ is the face poset of a CW sphere. If, in addition, $X$ is strongly regular (in other words, $X$ satisfies the intersection property), then $L(X)$ is a lattice. The face lattice of a $d$-simplex is called the {\bf Boolean lattice} of rank $d+1$, denoted by $B^{d+1}$. Every interval of a Boolean lattice is a Boolean lattice.

	Let $X$ be a strongly regular CW complex and let  $L=L(X)$ be its face lattice. For any $k$-face $F\in X$, the rank of $F$ in $L$ is $\rho(F) = k+1$. We denote by $f_k(X)$ or $f_k(L)$ the number of $k$-faces of $X$. In this paper, we will use the $f$-number notation for all graded lattices, not just those that are face lattices. (Note that $f_k(L) = W_{k+1}(L)$ where $W_k$ is the $k$-th Whitney number of the second type.) The {\bf pyramid} operation on a finite bounded poset $L$ (as defined in \cite{MR2335713}) is the Cartesian product with the Boolean lattice $B^1 = \{\hat{0}, \hat{1}\}$. That is, for a poset $L= [\hat{0}, \hat{1}]$ we let $\Pyr(L) := L \times B^1$ (with the componentwise partial order).

	Let $X$ be a regular CW complex and $L= L(X)$ be its face poset. The {\bf order complex} $O(X)$ (or $O(L)$) of $X$ is the following simplicial complex. The vertices of $O(X)$ are the nonempty proper faces of $X$, and the faces of $O(X)$ are the chains of nonempty proper faces of $X$. One advantage of working with a regular CW complex $X$ is that (the geometric realization of) the order complex $O(X)$ is homeomorphic to $X$ (see \cite{MR1373690}).

	A {\bf pseudomanifold} is a pure regular CW complex of dimension $d-1$ such that every $(d-2)$-face is contained in exactly two $(d-1)$-faces. It is easy to check that this property is preserved under taking barycentric subdivisions. A $(d-1)$-dimensional pseudomanifold $X$ is called {\bf normal} if $X$ is connected and $H_{d-1}(X)$ and $H_{d-1} (X, X - x)$ are $1$-dimensional for all $x \in X$. Here $H_i(X)$ and $H_i(X,X-x)$ are singular and relative homologies, respectively, computed with coefficients in in $\mathbb{Z}/2\mathbb{Z}$. For any simplicial complex $\Delta$ and a face $\tau\in \Delta$, the {\bf link of $\tau$ in $\Delta$}, denoted by $\lk_\Delta \tau$ is a subcomplex of $\Delta$ defined as follow.
	\begin{align*}
	\lk_\Delta \tau = \{\sigma \; |\; \sigma \cap\tau = \emptyset,\; \sigma \cup \tau \in \Delta   \}.
	\end{align*}
	If $X$ is a normal pseudomanifold and $T$ is a triangulation of $X$, then $\lk_T(\tau)$ is a normal pseudomanifold for each simplex $\tau \in T$.

	We have the following hierarchy relating face lattices of polytopes, strongly regular (SR) CW spheres/pseudomanifolds, and diamond lattices.
	\begin{equation}
	\{\text{polytopes} \}
	\subset \{\text{SR CW spheres} \}
	\subset \{
	\text{SR normal pseudomanifolds} \}
	\subset \{
	\text{diamond lattices} \}.\nonumber
	\end{equation}
	\normalsize


	\section{From polytopes to lattices}\label{section: extend lbt to lattices}
	Theorem 3.2 of \cite{myLBTpaper} established Gr\"unbaum's conjecture for all polytopes. The goal of this section is to extend this result to a much larger class of objects, namely, to all diamond lattices.
	
	First we point out that the key ingredient of the proofs in \cite{myLBTpaper} relied on the following property of polytopes: {\em if $P$ is a polytope and $u,v$ are two arbitrary distinct vertices of $P$, then there exists a facet $F$ such that $v\in F$, but $u\notin F$.}
	The notion of coatom-distinguishability introduced in Section 2 was motivated by this property of polytopes. 
	
	The following few propositions hold for all $\Diamond$-lattices. All of them can be easily proved by induction on rank. We include the proofs here for the sake of completeness but with fewer and fewer details as these proofs are very similar to each other.
	
	\begin{proposition}\label{prop: diamond implies coatom-dist} Every  $\Diamond$-lattice is coatom-distinguishable.
	\end{proposition}
	
	\begin{proof}
		Let $L = [\hat{0}, \hat{1}]$ be a $\Diamond$-lattice of rank $d+1$. We will prove the statement by induction on $d$. When $d=1$, $L$ is a diamond, and diamonds are coatom-distinguishable. Assume $d\geq 2$, and let $t,s\in L$ be elements of the same rank. If $s$ and $t$ are coatoms or if $s\vee t=\hat{1}$, there is nothing to prove. Otherwise, let $F$ be a coatom that is above $s \vee t$. Since the interval $[\hat{0},F]$ is a $\Diamond$-lattice of rank $d$, there must be a rank-$(d-1)$ element $G < F$ such that $t < G$ while $s\not< G$. Since the interval $[G,\hat{1}]$ has rank $2$, by the diamond property there exists another coatom $F'$ of $L$ such that $G = F \wedge F'$. This implies that $t < F'$ and $s\not< F'$.
	\end{proof}
	
	\begin{proposition}\label{prop: atleastBoolean}
		Let $L = [\hat{0}, \hat{1}]$ be a $\Diamond$-lattice of rank $d+1$. For every element $t\in L$, the interval $[t, \hat{1}]$ has at least the same number of elements in each rank as the Boolean lattice $B^{d+1-\rho(t)}$ (where $\rho(t)$ is the rank of $t$).
	\end{proposition}
	
	\begin{proof}
		The proof is again by induction on $d$ and $\rho(t)$. When $d=1,2$, the statement is clear. Let $d\geq 3$. If $\rho(t)=d+1$ or $d$ then $[t,\hat{1}] = B^0$ or $B^1$ respectively. If $\rho(t)=k<d$, there must exist (at least) two elements $s_1, s_2$ that cover $t$ in $[t,\hat{1}]$. Since $P$ is coatom-distinguishable by the previous proposition, there is a coatom $F$ such that $s_1 < F$ and $s_2\not< F$. By the inductive hypothesis, the two intervals $[t,F]$ and $[s_2, \hat{1}]$  both have at least the same number of elements in each rank as the Boolean lattice $B^{d-\rho(t)}$ and they are disjoint from each other. Since $f_{k}([t,\hat{1}]) \geq f_k([t,F]) + f_{k-1}([s_2,\hat{1}])$ for $k\geq 1$, the result follows by the Pascal triangle recursion. \end{proof}
	
	\begin{proposition}
		Let $L = [\hat{0}, \hat{1}]$ be a $\Diamond$-lattice of rank $d+1$. Consider any lower or upper rank-$k$ interval in $L$. If this interval has exactly $k$ atoms or exactly $k$ coatoms, then it is the Boolean lattice  $B^{k}$.
	\end{proposition}
	
	\begin{proof}
		We will, once again, induct on the rank of $L$. The base case is clear. Now let $\rho(L)=d+1>2$. Notice that if any proper lower or upper interval of $L$ satisfies the condition in the statement, then this interval is a Boolean lattice by the inductive hypothesis. Now suppose $L$ has $d+1$ atoms. By the lattice property of L, Each atom $x$ is below at most $d$ rank-$2$ elements since each rank-$2$ element has to be the join of $x$ and another atom. Together with Proposition \ref{prop: atleastBoolean}, this means each atom is below exactly $d$ rank-$2$ elements. By the inductive hypothesis, each proper upper interval of $L$ is a Boolean lattice. Now since $L$ has $d+1$ atoms, it is the Boolean lattice $B^{d+1}$. 
		
		The argument for the case when $L$ has $d+1$ coatoms is very similar, hence we omit it.
	\end{proof}
	
	\begin{proposition}\label{prop: all-cones are Boolean}
		Let $L = [\hat{0}, \hat{1}]$ be a $\Diamond$-lattice of rank $d+1$. Consider any lower or upper interval $[x,y]$ in $L$. If every coatom of this interval contains all the atoms of $[x,y]$ but one, then $[x,y]$ is a Boolean lattice.
	\end{proposition}
	
	\begin{proof}
		The proof is, again, by induction on the rank of the interval. The base case is trivial. Suppose the statement is true for all intervals of rank at most $k-1$ and let $[x,y]$ be a rank-$k$ interval satisfying the conditions in the statement above. Then for any coatom $z$ of $[x,y]$, the interval $[x,z]$ also satisfies conditions of the statement, and hence by the inductive hypothesis $[x,z] = B^{k-1}$. Since there is only one atom of $[x,y]$ that is not below $z$, this forces $[x,y]$ to have exactly $k$ atoms. By the previous proposition, $[x,y] = B^k$.
	\end{proof}

	Now we are ready to obtain the key ingredient of the proof of Gr\"unbaum's Conjecture, this time, for $\Diamond$-lattices. The following proofs are adaptations of the original proofs from \cite{myLBTpaper} to the generality of diamond lattices. We will only be using part (i) of Proposition \ref{thm: key prop} in this section. Parts (ii) and (iii) will become useful later in Section \ref{section: 2d+1 verts}.
	
	\begin{proposition}\label{thm: key prop}
		Let $L = [\hat{0},\hat{1}]$ be a $\Diamond$-lattice of rank $d+1$, and let $\mathcal{S}=\{v_1,v_2,\dots ,v_m\}$ be a subset of atoms of $L$. Let $|L_{k+1}^ \mathcal{S}|$ denote the number of rank-$(k+1)$ elements of $L$ that are above at least one of the $v_i$'s in $\mathcal{S}$. Then,
		\begin{itemize}
			\item[(i)] $|L_{k+1}^ \mathcal{S}| \geq \sum_{i=1}^{m} {d-i+1 \choose k}$.
			\item[(ii)] If one of the atoms in $\mathcal{S}$ is covered by more than $d$ elements of rank $2$, then $|L_{k+1}^ \mathcal{S}| \geq \bigg( \sum_{i=1}^{m} {d+1 -i \choose k} \bigg) + {d-2 \choose k-1}$.
			\item[(iii)] If there exists a coatom $F\in L$ such that $F$ is above at least one and not more than $|S|-2$ elements of $S$, then  $|L_{k+1}^ \mathcal{S}| \geq  {d \choose k} + {d-1 \choose k}+ \sum_{i=1}^{m-2} {d -i \choose k}$.
		\end{itemize}
	\end{proposition}

	\begin{proof}
		We will only prove part (i) since the proofs of (ii) and (iii) are essentially the same.
		
		We induct on $m$ to show that there exists a sequence of elements, $\{F_1,\;\dots , \; F_m \}$, such that 
		\begin{enumerate}
			\item[(1)] each $F_i$ has rank $d-i+2$, 
			\item[(2)] $v_i \leq F_i$, but $v_j \not\leq F_i$ for any $j<i$.
		\end{enumerate}
		
		The base case is $m=1$, and we simply pick $F_1 = \hat{1}$. Inductively we assume that for every $p\leq m-1$ and any $p$-element set of vertices $\{v_1, \;\dots,\; v_p\}$, there exists a sequence $\{F_1, \; \dots ,\; F_p\}$ such that for $1\leq i\leq p$, conditions (1) and (2) are satisfied.
		
		Let $m>1$ and let $v_1,\; \dots , \; v_m$ be $m$ given atoms of $L$. By the inductive hypothesis, for $\{v_1,\;\dots ,\;v_{m-1} \}$ there exist elements $F_1,\;\dots ,\; F_{m-1}$ satisfying conditions (1) and (2). Similarly, by considering $\{v_1,\;\dots ,\;v_{m-2}, v_m \}$, there also exists a rank-$(d-m+3)$ element $F$ such that $v_m \leq F$ but $v_1,\;\dots,\; v_{m-2} \not\leq F$. Regardless of whether $v_{m-1}$ is below $F$ or not, by the coatom-distinguishability of $L$, there must exist a coatom of $[\hat{0},F]$, call it $F_m$, that satisfies $v_m\leq F_m$ but $v_{m-1} \not\leq F_m$. Then $v_i\in F_m$ if and only if $i = m$, and $F_1,\;\ldots,\; F_{m-1}, F_m$ is a desired sequence.
		
		For each $i$, the rank-$(k+1)$ elements in $F_i$ that are above $v_i$ correspond to rank-$k$ elements in the interval $[v_i, F_i]$. From this we obtain that
		\begin{align*}
		&\quad\; \#\; \text{rank-}(k+1) \text{ elements of } L \text{ that contain some } v_i \; (1\leq i\leq m) \\
		&\geq 	 \# \;\bigcup_{i=1}^{m} \{ \text{rank-}(k+1)\text{ elements of }  F_i \text{ containing } v_i \}\\
		&= \#\; 	\bigcup_{i=1}^{m} \{ \text{rank-}k\text{ elements} \text{ in the lattice } [v_i, F_i]\}\\
		&\geq  \sum_{i=1}^{m}{d-i+1 \choose k},
		\end{align*} 
		where the last inequality follows from Proposition \ref{prop: atleastBoolean}. The result follows.
	\end{proof}
	
	With Proposition \ref{thm: key prop} at our disposal, we are now ready to prove the first main theorem of this section. Recall that 
	\[\phi_k(d+s,d) = {d+1 \choose k+1 }+{d \choose k+1 }-{d+1-s \choose k+1 }.\]
	
	\begin{theorem}\label{thm: inequality lattice}
		Let $L$ be a rank-$(d+1)$ $\Diamond$-lattice that has $d+s$ atoms where $d\geq s\geq 1$. Then $f_k(L)\geq \phi_k(d+s,d)$ for every $k$.
	\end{theorem}
	
	\begin{proof}
		The statement holds for $s=1$ by Proposition \ref{prop: atleastBoolean}. The proof is by induction on $s$. We fix $s\geq 2$. The following argument will show that if the statement holds for all pairs $(s', d')$ such that $s'<s$ and $d'\geq s'$, then it also holds for the pairs $(s,d)$ for all $d\geq s$. 
		Now consider $d\geq s$, and let $L=[\hat{0},\hat{1}]$ be a rank-$(d+1)$ $\Diamond$-lattice with $d+s$ atoms. If there exists a coatom $z$ of $L$ such that  $d\leq f_0([\hat{0},z])\leq d+s-2$, or equivalently, if $[\hat{0},z]$ is a $\Diamond$-lattice with $(d+s-m)$ atoms such that $2\leq m\leq s$, then there are $m$ atoms of $L$ that are not below the element $z$. We denote $m$ of them by $\{v_1, v_2, \dots , v_m \}$. The rank-$k$ elements of $L$ fall into two disjoint categories: those that are below $z$ and those that are above some $v_i$. By the inductive hypothesis on $[\hat{0}, z]$, $f_k([\hat{0}, z]) \geq \phi_k(d+s-m,d-1)$. Therefore by Proposition \ref{thm: key prop},
		\begin{eqnarray}\label{eq: d+s proof}
		f_k(P) &\overset{(\star)}{\geq}& \phi_k(d+s-m,d-1) + \sum_{i=1}^{m} {d-i+1 \choose k}\nonumber 
		\\
		&=& \phi_k(d+s-m+2,d) + \sum_{i=3}^{m}  {d-i+1 \choose k} \nonumber 
		\\
		&=& \left[\phi_k(d+s,d) - \sum_{j=1}^{m-2} {d-s+m-1-j \choose k} \right] + \sum_{i=3}^{m} {d-i+1 \choose k}
		\\
		&=& \phi_k(d+s,d) + \sum_{j=1}^{m-2} \underbrace{\left[ -{d-j-1-(s-m) \choose k}+ {d-j-1 \choose k}\right]}_{\geq 0 \text{ (since } s\geq m)} \nonumber 
		\\ 
		&\overset{(\star\star)}{\geq}& \phi_k(d+s,d).\nonumber
		\end{eqnarray}
		This completes the proof of this case.
		
		Otherwise, all coatoms in $P$ contain all but one atom. By Proposition \ref{prop: all-cones are Boolean} the interval $[x,y]$ is the Boolean lattice.
	\end{proof}

	Using the inequalities $(\star)$ and $(\star\star)$ in the proof above, we can obtain the following corollary. The proof is the same as in \cite{myLBTpaper} hence we omit it.
	
	\begin{corollary}\label{cor:diamonds}
		If $f_k(L)=\phi_k(d+s,d)$ for some $k$ with $1\leq k\leq d-2$, then each coatom of $L$ lies above $d$, $d+s-2$, or $d+s-1$ atoms, and $L$ has $d+2$ coatoms.
	\end{corollary}
	
	Our next goal is to extend the treatment of equality cases from the class of polytopes to the class of strongly regular normal pseudomanifolds. The next main theorem in this section is the following. We denote by $\Pyr^{d-s}(\Delta^1 \times \Delta^{s-1})$ the polytope that is a $(d-s)$-fold pyramid over a prism over an $(s-1)$-simplex.
	
	\begin{theorem} \label{thm:equality sphere}
		Let $L$ be the face lattice of a strongly regular normal $(d-1)$-pseudomanifold $S$ with $d+s$ vertices where $s\geq 2$ and $d\geq s$. If $f_k(L)=\phi_k(d+s,d)$ for some $1\leq k\leq d-2$, then $S$ is the boundary of $\Pyr^{d-s}(\Delta^1 \times \Delta^{s-1})$.
	\end{theorem}
	
	By \cite[Thm.~4.3]{myLBTpaper}, the statement holds when $S$ is (the boundary of) a polytope. By Corollary \ref{cor:diamonds}, $S$ has $d+2$ facets. Therefore it suffices to show that all strongly regular normal $(d-1)$-pseudomanifolds with $d+2$ facets are boundary complexes of polytopes. Before proving this, we will prove the following useful proposition.
	
	\begin{proposition}\label{prop: simplicial dual is a sphere}
		Let $L=L(P)$ be the face lattice of a regular CW complex $P$. If all proper upper intervals of $L$ are Boolean lattices, then its dual lattice $L^*$ is the face lattice of a simplicial complex that is homeomorphic to $P$ .
	\end{proposition}
	
	\begin{proof}
		Since all proper upper intervals of $L$ are Boolean lattices, so are all proper lower intervals of $L^*$. Therefore $L^*$ is the face lattice of a pure, connected simplicial complex $\Delta$. Now the two order complexes $O(L)$ and $O(L^*)$ are homeomorphic to $P$ and $\Delta$ respectively. Since $O(L)=O(L^*)$, it follows that $P$ and $\Delta$ are homeomorphic.
	\end{proof}
	
	We can now prove the following theorem, which is the last piece in the proof of Theorem \ref{thm:equality sphere}.
	
	\begin{theorem}\label{thm: spheres are polytopes}
		Every strongly regular $(d-1)$-dimensional normal pseudomanifold that has $d+2$ facets is combinatorially equivalent to the boundary complex of a $d$-polytope.
	\end{theorem}

	\begin{proof} 
		
		Let $P$ be a strongly regular $(d-1)$-dimensional normal pseudomanifold that has $d+2$ facets. This implies that every $k$-face of $P$ has at most $k+2$ facets. We will, one more time, induct on $d$. The base cases ($d = 2$) are trivial. Now for $d>2$, we first consider the case when $L(P)$ is a pyramid over $L(F)$ for some facet $F$ of $P$. By the inductive hypothesis $F$ is combinatorially equivalent to a polytope, hence so is $P$.
		
		Next we assume that $P$ is not a pyramid, which implies that, for any vertex $v$ there exist at least two facets that do not contain $v$. Let $L$ be the face lattice of $P$. We now look at its dual lattice, $L^*$, which has $d+2$ atoms. Since $P$ is not a pyramid, for every coatom $v\in L^*$, the lower interval $[\hat{0}^*, v]$ has rank $d$ and at most $d$ atoms. This forces $[\hat{0}^*, v]$ to be the Boolean lattice, and by Proposition \ref{prop: simplicial dual is a sphere}, $L^*$ is the face lattice of a simplicial normal $(d-1)$-pseudomanifold with $d+2$ vertices, which must be the join of the boundary complexes of two simplices whose dimensions add up to $d$ (see, for example, \cite[Lemma 4.1]{Hailunspaper}). Thus $L^*$ is the face lattice of a polytope, and hence so is $L$.
		
	\end{proof}

	\section{Polytopes with $d+2$ facets}\label{section: d+2 facets}
	
	 In the previous section, we extended all the results from \cite{myLBTpaper} about polytopes with up to $2d$ vertices to the generality of diamond lattices and strongly regular normal pseudomanifolds. Our goal in the rest of the paper is to provide another generalization: this time to strongly regular normal pseudomanifolds with $2d+1$ vertices. To do so, in this section we will temporarily go back to the world of polytopes. We start by using the face numbers to define a certain partial order on the $d$-polytopes.
	
	Given two $d$-polytopes $P$ and $Q$, we write $P \leq Q $ if $f_i(P) \leq f_i(Q)$ for all $0\leq i\leq d-1$; if, in addition, $f(P) \neq f(Q)$, then we write $P<Q$. If $\mathcal{S}$ is a set of $d$-polytopes, and there exists a polytope $P\in \mathcal{S}$ such that for all other polytopes $Q\in \mathcal{S}$, $P<Q$, then $P$ is called {\bf componentwise minimal} in $\mathcal{S}$ and $f(P)$ is called the {\bf componentwise minimal $f$-vector in $\mathcal{S}$}. If moreover we are able to order all polytopes in $\mathcal{S}$ so that $\mathcal{S} =\{P_1 < P_2 < P_3 < \dots \}$, then $\mathcal{S}$ is {\bf completely ordered by $f$-vectors}. Notice that this is a very strong condition, and it implies that any nonempty subset of $\mathcal{S}$ contains a polytope that is componentwise minimal.
	
	Let $\mathscr{P}(d+s,d, d+2)$ be the set of all  $d$-polytopes that have $d+s$ vertices and $d+2$ facets. The main result in this section is the following.
	\begin{theorem}\label{thm: remark}
		Polytopes in $\mathscr{P}(d+s,d, d+2)$ are completely ordered by $f$-vectors.
	\end{theorem}
	
	The set $\mathscr{P}(d+s,d, d+2)$ consists of polytopes dual to those with $d+2$ vertices and $d+s$ facets. The set of all $d$-polytopes with $d+2$ vertices, denoted $\mathscr{P}(d+2,d)$, was fully characterized in \cite[6.1]{Gr1-2}. In particular, every simplicial polytope in this class is obtained from some $(d-m)$-simplex by adding one additional vertex and placing it beyond exactly $m$ of its facets for some $1\leq m\leq \lfloor\frac{d}{2}\rfloor$. Gr\"unbaum denoted such a simplicial polytope as $T^d_m$. Moreover, every $P\in \mathscr{P}(d+2,d)$ is a $(d-i)$-fold pyramid over $T^i_m$ for some pairs $(i,m)\in S := \{2\leq i\leq d,\; 1\leq m\leq \lfloor\frac{i}{2}\rfloor\}$. Such a polytope is denoted as $T^{d,d-i}_m$. Gr\"unbaum showed in  \cite[6.1.4]{Gr1-2} that,
	\begin{align}\label{eqn: f_k}
	f_k(T^{d,d-i}_m) ={d+2 \choose d-k+1} - {d-i+m+1 \choose d-k+1} - {d-m+1 \choose d-k+1} + {d-i+1 \choose d-k+1} \end{align} for all $0\leq k\leq d-1$. In particular, 
	\begin{align}\label{eqn: f_d-1}
	f_{d-1}(T^{d,d-i}_m) = d + 1 + m(i-m).\end{align}
	
	Equation (\ref{eqn: f_k}) implies the following inequalities on the face numbers of polytopes in $\mathscr{P}(d+2,d)$ (also see \cite[6.1]{Gr1-2}): \begin{itemize}
		\item[(i)] For $(i,m)$, $(i,m+1) \in S$, 
		\begin{align}\label{relation 1} f_k(T^{d,d-i}_m) \leq f_k(T^{d,d-i}_{m+1}), \end{align}
		with strict inequality if and only if $m\leq k+1$.
		\item[(ii)] For $(i,m) \in S$, 
		\begin{align}\label{relation 2} f_k(T^{d,d-i}_m) \leq f_k(T^{d,d-i-1}_{m}), \end{align}
		with strict inequality if and only if $d\leq d-i+ k+m$.
	\end{itemize}
	
We are now ready to prove Theorem \ref{thm: remark}.
	\begin{proof}[Proof of Theorem \ref{thm: remark}]
		Let $\leq$ be the partial order on polytopes defined at the start of this section. Note that if two polytopes $P$ and $Q$ satisfy $P\leq Q$, then so do their duals, that is, $P^*\leq Q^*$. Thus to prove the statement, it suffices to show that $\mathscr{P}(d+2,d,d+s)$ is completely ordered by $f$-vectors. By (\ref{eqn: f_d-1}), for any $i\neq i'$, $T^{d,d-i}_m $ and $T^{d,d-i'}_m$ have different numbers of facets. Therefore two (distinct) polytopes $P = T^{d,d-i}_m$ and $P'= T^{d,d-i'}_{m'}$ are both in $\mathscr{P}(d+2,d, d+s)$ only if $m\neq m'$ (similarly, if $m=m'$, then $i= i'$ for polytopes in this set). This implies that each polytope $P$ in $\mathscr{P}(d+2,d,d+s)$ is of the form $T_m^{d,d-i}$ where all $m$ (and all $i$) are distinct (that is, the functions $P\to m$ and $P\to i$ are one-to-one). Hence, in order to prove the statement, it suffices to show that if $\mathscr{P}(d+2,d, d+s)$ contains both $T^{d,d-i'}_{m'}$ and $T^{d,d-i}_{m}$ with $m < m'$, then
		\[f_k(T^{d,d-i}_{m})\leq f_k(T^{d,d-i'}_{m'})\quad \text{ for any } 1\leq k\leq d-2.\]
		
		Since $T^{d,d-i}_m$ and $T^{d,d-i'}_{m'}$ have the same number of facets, and the inequalities in (\ref{relation 1}) and (\ref{relation 2}) are both strict when $k=d-1$, the fact that $m<m'$ implies $i>i'$.
		
		For any $1\leq k\leq d-2$, by (\ref{eqn: f_k}),
		\begin{align*}
		&f_{k}(T^{d,d-i'}_{m'}) - f_{k}(T^{d,d-i}_{m})\\
		&= \underbrace{\left[- {d-i'+m'+1 \choose d-k+1} +{d-i+m+1 \choose d-k+1}\right]}_{(=A)}  + \underbrace{\left[ - {d-m'+1 \choose d-k+1} + {d-m+1 \choose d-k+1} \right]}_{(=B)} \\
		&\quad\;\;  + \underbrace{\left[{d-i'+1 \choose d-k+1} - {d-i+1 \choose d-k+1}\right]}_{(=C)}.
		\end{align*}
		By repeatedly applying the Pascal identity, we obtain the following (known) formula.
		\[{n \choose k+1} - {n-a \choose k+1} = \sum_{p=1}^{a} {n-p \choose k}.\] Therefore
		\begin{align*}A = - \sum_{p=1}^{(i-i')+(m'-m)} {d-i'+m'+1-p \choose d-k}, 	\end{align*}
		
		\begin{align*}B =\sum_{p=1}^{m'-m} {d-m+1-p \choose d-k},
		\end{align*}
		\begin{align*}
		C = \sum_{p=1}^{i-i'} {d-i'+1-p \choose d-k}.
		\end{align*}
		By shifting the summation index of $A$ and then cancelling with the overlapping terms in $C$, we obtain
		\begin{align*}\label{eqn: A+B+C}
		A + C &= - \sum_{p=1-m'}^{0} {d-i'+1-p \choose d-k} + \sum_{p=(i-i')-m+1}^{i-i'} {d-i'+1-p \choose d-k}\nonumber \\
		&= -\sum_{p=1}^{m'} {d-i'+p \choose d-k} + \sum_{p=1}^{m} {d-i+p \choose d-k}\nonumber \\
		&= -\sum_{p=1}^{m'-m} {d-i'+p+m \choose d-k} - \sum_{p=1}^{m} \left[{d-i'+p \choose d-k} - {d-i+p \choose d-k}  \right],\text{ while} \nonumber 
		\end{align*}
		\begin{align}
		B+A+C &= \sum_{p=1}^{m'-m} \left[{d-m+1-p \choose d-k} -{d-i'+p+m \choose d-k} \right] - \sum_{p=1}^{m} \left[{d-i'+p \choose d-k} - {d-i+p \choose d-k}  \right] \nonumber\\
		&= \sum_{p=1}^{m'-m}\left[\sum_{q=1}^{i'+1-2p-2m}{d-m+1-p-q \choose d-k-1}\right] - \sum_{p=1}^{m} \left[\sum_{q=1}^{i-i'}{d-i'+p-q \choose d-k-1}\right] \\
		&\overset{(\bowtie)}{\geq }\left[\sum_{p=1}^{m'-m}(i'+1-2p-2m) \right]{d-m+1-(1) -(i'-1-2m) \choose d-k-1} \nonumber \\
		&\qquad -(m)(i-i'){d-i'+m-1 \choose d-k-1} \nonumber \\
		&= (m'-m)(i'-m'-m){d-i'+m+1 \choose d-k-1} -(m)(i-i'){d-i'+m-1 \choose d-k-1}. \nonumber
		\end{align}
		By (\ref{eqn: f_d-1}), \[d+s = d+1+ m(i-m) = d+1+ m'(i'-m'),\]
		which implies that $m(i-m) = m'(i'-m')$, and so 
		\begin{align}\label{eqn: (m, i) vs (m', i') }
		(m'-m)(i'-m'-m)=(m)(i-i').
		\end{align}
		Therefore
		\begin{align}\label{eqn: end}
		f_k(T^{d,d-i'}_{m'})-f_k(T^{d,d-i}_{m}) =B+A+C &\geq \underbrace{(m)(i-i')}_{>0}\underbrace{\left[{d-i'+m+1 \choose d-k-1} -{d-i'+m-1 \choose d-k-1} \right]}_{\geq 0 }\\
		&\geq 0.\nonumber
		\end{align}
		The result follows.
	\end{proof}
	
	 Recall that the simplicial polytope $T^{s}_{1}$ is a bipyramid over an $(s-1)$-simplex, and $T^{d,d-s}_1$ is a $(d-s)$-fold pyramid over this bipyramid. When $s\leq d$, the $d$-polytope $T^{d,d-s}_1$ exists and is contained in the set $\mathscr{P}(d+2,d,d+s)$. By Theorem \ref{thm: remark}, it has the componentwise minimal $f$-vector in the completely ordered set $\mathscr{P}(d+2,d,d+s)$. 
Together with the fact that $f_k\left((T^{d,d-s}_1)^*\right) = \phi_k(d+s,d)$, this implies the following corollary.
	
	\begin{corollary}
		For any $P\in \mathscr{P}(d+s,d, d+2)$, $2\leq s\leq d$, and $1\leq k\leq d-2$, \[f_k(P) \geq \phi_k(d+s,d).\] 
	\end{corollary}

	When $s>d$, the polytope $T^{d,d-s}_1$ does not exist. By Theorem \ref{thm: remark} the ``next-best'' polytopes would be $T^{d,d-i}_2$ (for some $i$) when they exist. If $d$ is even and $s=d+1$, then the polytope $T^{d,\frac{d}{2}-2}_2$ has the componentwise minimal $f$-vector in $P(d+2, d, 2d+1)$. However, when $d$ is odd, $\frac{d}{2}$ is not an integer. Yet in the proof above, the inequality (\ref{eqn: end}) continues to hold even though the equality in (\ref{eqn: (m, i) vs (m', i') }) becomes ``$\geq$''. Hence we only need to slightly adjust the $f$-vector formula of  $T^{d,\frac{d}{2}-2}_2$ by using the floor or ceiling function as needed to keep the inequalities in the proof valid, so the proof above still works, and the face numbers of hypothetical $T^{d,\frac{d}{2}-2}_2$ (obtained from (\ref{eqn: f_k}) and shown below) form the lower bounds.
	
	\begin{lemma}\label{lem: bound for (2d+1, d+2)}
		Let $P$ be a $d$-polytope with $f_0(P)\geq 2d+1$ and $f_{d-1}(P) = d+2$. Then for $m\geq 1$, 
		\begin{align}\label{eq: bound for (2d+1, d+2)}
		f_m(P) \geq {d+1 \choose m+1} + {d \choose m+1} + {d-1 \choose m+1} - {\ceil{\frac{d}{2}} \choose m+1} - {\ceil{\frac{d}{2}} -1 \choose m+1}.
		\end{align}
		When $f_0(P) = 2d+1$ and  $d$ is even, this lower bound is sharp, and it is attained for every $m\geq 1$ when $P  = (T^{d,\frac{d}{2}-2}_2)^*$. 
	\end{lemma}


	\section{Pseudomanifolds with at least $d+3$ facets}\label{section: 2d+1 verts}
	In this section, we will find the lower bounds on the $f$-numbers of polytopes and strongly regular normal pseudomanifolds with at least $2d+1$ vertices and $d+3$ facets.  In particular, we will show that in this case, there is a polytope $\nabla$ (constructed below) that has the componentwise minimal $f$-vector. 
	
	Let $P$ be a polytope that has a facet $F$ that is a simplex. The operation that builds a shallow pyramid over $F$ while maintaining the convexity is called {\bf stacking} on $P$ over $F$. Let $T^{d,d-2}_1$ be the polytope with $d+2$ vertices as defined in Section \ref{section: d+2 facets}. Stack over one arbitray facet of $T^{d,d-2}_1$ that is a simplex (note that $T^{d,d-2}_1$ is a $(d-2)$-fold pyramid over a square), then take the dual of this polytope. We will obtain the polytope $\nabla := \big(\Stack(T^{d,d-2}_1)\big)^*$, which has $d+3$ facets. It is also easy to check that $f_0(\nabla)=2d+1$ and
	\begin{align}
	f_m(\nabla) = {d+1 \choose m+1} + {d \choose m+1} + {d-1 \choose m} \quad \text{for } 0<m \leq d-1.
	\end{align}
	The main theorem of this section is the following.
	\begin{theorem}\label{thm: 2d+1}
		Let $P$ be a $(d-1)$-dimensional strongly regular normal pseudomanifold with $f_0(P)\geq 2d+1$ and $f_{d-1}(P)\geq d+3$. Then $P \geq \nabla$.
	\end{theorem}
	
The proof of this theorem consists of several cases. We first use the Lower Bound Theorem to take care of the case when $P$ is simple, i.e., every vertex is in $d$ edges.
	\begin{theorem}\label{thm:2d+1, d+3, simple case}
		Let $P$ be a pseudomanifold as above. If, in addition, $P$ is simple, then $P \geq \nabla$.
	\end{theorem}
	\begin{proof}
		Let $L=L(P)$ be the face lattice of $P$. Since every vertex of $P$ is in exactly $d$ edges, the dual lattice $L^*$ is the face lattice of a pure simplicial complex. Since $L$ is the face lattice of a normal pseudomanifold, by Proposition \ref{prop: simplicial dual is a sphere}, $L^*$ is the face lattice of a simplicial normal pseudomanifold, call it $P^*$. Since $f_{0}(P^*) \geq d+3$, by the Lower Bound Theorem (see \cite{MR2636638}, \cite{MR877009}, and \cite{MR1314963}), the $f$-vector of $P^*$ is componentwise bounded by the $f$-vector of the stacked polytope with $d+3$ vertices, denoted by $\Stack(d+3,d)$. For each $m\geq 1$,
		\begingroup
		\allowdisplaybreaks
		\begin{align}\label{ineq: case 1}
		f_m(P) &\geq f_{d-1-m}\big(\Stack(d+3,d)\big)\nonumber\\
		& = {d \choose m+1 }\cdot (d+3) - {d+1 \choose m+1 }\cdot (d-1-m) \nonumber\\
		& = {d+1 \choose m+1 } + {d \choose m+1 } +\bigg[{d \choose m+1 }\cdot (d+2) - {d+1 \choose m+1 }\cdot (d-m) \bigg]   \nonumber\\
		& =  {d+1 \choose m+1 } + {d \choose m+1 } + {d \choose m+1} \\
		& = f_m(\nabla) + {d-1 \choose m+1} \nonumber\\
		&\geq  f_m(\nabla). \nonumber
		\end{align}
		\endgroup
		The result follows.
	\end{proof}

	\begin{proof}[Proof of Theorem \ref{thm: 2d+1}] 
		Let $P$ be a $(d-1)$-dimensional strongly regular normal pseudomanifold with $f_0(P)= 2d+k$ ($k\geq 1$) and $f_{d-1}(P)\geq d+3$. By Theorem \ref{thm:2d+1, d+3, simple case}, we now only need to consider the case where there is a vertex $v\in P$ that is {\bf nonsimple}, i.e., $v$ is in at least $d+1$ edges. In the rest of the proof $v$ always denotes a non-simple vertex of $P$. The proof that Theorem \ref{thm: 2d+1} holds in this case is by induction on $d$.
		
		The statement clearly holds for $d=2$. For the inductive step, we treat the following four cases.
		\begin{itemize}
			\item[Case 1:] $P$ is the pyramid with apex $v$, that is, there exists a facet $F_1$ of $P$ such that $P = F_1 * v$.
			In this case, the face lattice of $F_1$ is the face lattice of a strongly regular CW $(d-2)$-sphere with at least $d+2$ facets and $2d+k-1 (> 2(d-1)+1)$ vertices. By the inductive hypothesis, for $m\geq 2$, we have
			\begin{align}
			f_m(P) &= f_m(F_1) + f_{m-1}(F_1) \nonumber\\
			&\geq \bigg[{d \choose m+1} + {d-1 \choose m+1} + {d-2 \choose m}    \bigg] + \bigg[{d \choose m} + {d-1 \choose m} + {d-2 \choose m-1}    \bigg]\nonumber\\
			&= {d+1 \choose m+1} + {d \choose m+1} + {d-1 \choose m}  = f_m(\nabla).\nonumber
			\end{align}
			
			And for $m=1$, 
			\begin{align}
			f_1(P) &= f_1(F_1) + f_{0}(F_1)\nonumber\\
			&\geq \bigg[{d \choose 2} + {d-1 \choose 2} + {d-2 \choose 1}    \bigg] + 2d+k-1 \nonumber\\
			&= {d+1 \choose 2} + {d \choose 2} + {d-1 \choose 1}  + \underbrace{k-1}_{\geq 0} \nonumber\\
			&\geq f_1(\nabla).\nonumber
			\end{align}

			\item[Case 2:] There exists a facet $F_1$ of $P$ that does not contain $v$ and such that $2d-1\leq f_0(F_1)\leq f_0(P)-2$ and $f_{d-2}(F_1)\geq (d-1)+3$. Then $F_1$ satisfies the inductive hypothesis. Furthermore, there are at least two vertices of $P$ outside of $F_1$ including a non-simple vertex $v$. By Proposition \ref{thm: key prop}(ii) we have that for all $m\geq 1$,
			\begin{align}\label{eqn: case 2}
			f_m(P) &\geq f_{m}(F_1) + {d \choose m} + {d-1 \choose m}+ {d-2 \choose m-1}\nonumber\\
			& \geq \underbrace{\bigg[ {d \choose m+1} + {d-1 \choose m+1}+ {d-2 \choose m} \bigg ]}_{(*)}+ {d \choose m}  + {d-1 \choose m}+ {d-2 \choose m-1}\\
			& = f_m(\nabla).\nonumber
			\end{align}
			
			\item[Case 3:] There exists a facet $F_1$ of $P$ that does not contain $v$ and such that $2d-1\leq f_0(F_1)\leq f_0(P)-2$ and $f_{d-2}((F_1) = (d-1)+2$. In this case, by Theorem \ref{thm: spheres are polytopes}, $F_1$ is combinatorially equivalent to (the boundary of) a polytope. By Lemma \ref{lem: bound for (2d+1, d+2)}, for $m\geq 1$, 
			\begin{align}\label{eq: bound for (2d+1, d+2)}
			f_m(F_1) \geq {d \choose m+1} + {d-1 \choose m+1} + {d-2 \choose m+1} - {\ceil{\frac{d-1}{2}} \choose m+1} - {\ceil{\frac{d-1}{2}} -1 \choose m+1}.
			\end{align}
			
			Recall that as in the previous case, there are two vertices outside of $F_1$ and one of them is nonsimple. Thus, using (\ref{eqn: case 2}), when $1\leq m< \ceil{\frac{d}{2}} -1 $ it suffices to show that $f_m(F_1)\geq {d \choose m+1} + {d-1 \choose m+1} + {d-2 \choose m}$. Recall that any polytope with $n$ vertices and $d+2$ facets is of the form $(T^{d,d-i}_m)^*$ for some $2\leq i\leq d$, $m<\floor{\frac{i}{2}}$, and $n = d+1+ m(i-m)$. It is easy to see that there does not exist a $5$-dimensional polytope with $11$ vertices and $7$ facets. Therefore when $d=6$, the facet $F_1$ has at least $12$ vertices. According to the ordering of all $5$-polytopes with $7$ facets as in Theorem \ref{thm: remark}, the minimizer should be $(T^{5,0}_2)^*$, hence $f_m(F_1)\geq f_m((T^{5,0}_2)^*)$. In this case if $m< \ceil{\frac{5}{2}} -1  =2$, then $m=1$. The result then follows since
			
			\begin{align*}
				{d \choose m+1} + {d-1 \choose m+1} + {d-2 \choose m} = 29 <30 = f_1((T^{5,0}_2)^*) \leq f_1(F_1).
			\end{align*}
			
			 When $d\geq 7$ and $m< \ceil{\frac{d-1}{2} }-1$, the lower bound in (\ref{eq: bound for (2d+1, d+2)}) is greater than ${d \choose m+1} + {d-1 \choose m+1} + {d-2 \choose m}$ in (\ref{eqn: case 2}). The computations to verify this are purely combinatorial, hence included in Appendix. 
			 
			 Now it suffices to focus on the cases when $m\geq  \ceil{\frac{d-1}{2} }-1$. If $f_0(F_1) \leq f_0(P)-3$, then there exist at least three vertices outside $F$. By Proposition \ref{thm: key prop},
			\begin{align*}
			f_m(P) &\geq f_m(F_1) + {d \choose m}  + {d-1 \choose m} + {d-2 \choose m} + {d-2 \choose m-1}\\
			&\geq  \bigg[ {d \choose m+1} + {d-1 \choose m+1} + {d-2 \choose m+1}  -  {\ceil{\frac{d-1}{2} }-1 \choose m}\bigg]+ {d \choose m}  + 2{d-1 \choose m}\\
			&=  {d+1 \choose m+1} + {d \choose m+1} + {d-1 \choose m}  + \underbrace{\bigg[ {d-2 \choose m+1}   - {\ceil{\frac{d-1}{2} }-1 \choose m}\bigg] }_{\geq 0 \text{ since } m \geq \ceil{\frac{d-1}{2} }-1}\\
			&\geq f_m(\nabla).
			\end{align*}
			
			Assume now that $f_0(F_1)=f_0(P)-2$, and let $v$ and $w$ be the only two vertices of $P$ outside of $F_1$. If $v$ has at least $d+2$ neighbors, then the vertex figure $P\slash v$ has at least $d+2$ vertices. Then we use similar arguments and computations in Proposition \ref{thm: key prop} to obtain:
			\[\#\{m\text{-faces of } P \text{ containing } v\} \geq {d \choose m } + {d-1 \choose m }  - {d-3 \choose m}. \]
			Therefore 
			\begin{align*}
			f_m(P) &\geq f_m(F_1) + {d \choose m}  + 2{d-1 \choose m}- {d-3 \choose m} \\
			&\geq  \bigg[ {d \choose m+1} + {d-1 \choose m+1} + {d-2 \choose m+1} - {\ceil{\frac{d-1}{2} }-1 \choose m} \bigg]+ {d \choose m}  + 2{d-1 \choose m} - {d-3 \choose m}\\
			&=  {d+1 \choose m+1} + {d \choose m+1} + {d-1 \choose m}  + \bigg[ {d-2 \choose m+1} - {d-3 \choose m} -  {\ceil{\frac{d-1}{2} }-1 \choose m}\bigg]\\
			&=  {d+1 \choose m+1} + {d \choose m+1} + {d-1 \choose m}  + \underbrace{\bigg[ {d-3 \choose m+1} -  {\ceil{\frac{d-1}{2} }-1 \choose m}\bigg] }_{\geq 0 \text{ since } m\geq \ceil{\frac{d-1}{2} }-1}\\
			&\geq f_m(\nabla).
			\end{align*}
			Otherwise, $v$ and $w$ each have at most $d+1$ neighbors. We will now show that this is impossible. First notice that $\{v,w\}$ must be an edge. 
			Since $f_0(F_1) = f_0(P) -2 \geq 2d-1$, there exist at least $2d-1$ edges connecting vertices of $F_1$ to $v$ and $w$. If at least $ d+1$ of these edges go to either $v$ or $w$, then since $\{v,w\}$ is also an edge, one of $v,w$ has at least $d+2$ neighbors. This takes us back to the previous case. Therefore there has to be at most $2d$ edges going out of $F_1$ (hence $f_0(F_1)=2d-1$ or $2d$). Since $v$ has degree $d+1$, $d$ of the vertices of $F_1$ will have edges connecting to $v$, and the other vertices of $F$ must be connected to $w$. Therefore there exists at most one vertex of $F_1$ that is connected to both $v$ and $w$. On the other hand, since every other facet of $P$ must intersect with $F_1$, but $f_{d-2}(F_1)=d+1$ while $f_{d-1}(P)\geq d+3$, there has to be at least one facet $F_2$ that does not intersect with $F_1$ at a ridge. The only possibility is that $F_2\cap F_1 = G$ is a $(d-3)$-dimensional, and $F_2 = G *v * w$. This means that every vertex of $G$ has an edge to both $v$ and $w$. By our assumptions it must be $f_0(G)=1$, and so $d=3$, $f_0(P)=7$ and $F_1$ is a pentagon.  Recall that $f(\nabla)= (7,11,6)$ when $d=3$. Any $3$-polytope with $7$ vertices that has a pentagon ($F_1$) and a triangle ($F_2$) that intersect at a vertex will need at least $12$ edges and seven $2$-faces. Therefore $P>\nabla$. 
			
			We are now done with Case 3. From now on we can assume that {\bf for each nonsimple vertex $v\in P$, each facet that does not contain $v$ has size at most $2d-2$.}
			
			\item[Case 4: ] There is a facet $F_1$ such that $v\notin F_1$ and $d+1 \leq f_0(F_1) = (d-1)+s \leq 2d-2$, i.e., $2\leq s\leq d-1$. First label the vertices of $P$ that are not in $F_1$ as $v= x_1, x_2, \dots, x_{d-s+2}$. By Proposition \ref{thm: key prop},
			\begin{align}\label{ineq: case 5 new}
			f_m(P)  &= f_m(F) + \#\{m \text{-faces of } P \text{ that contain some } x_i \}\nonumber\\
			&= \bigg[{d \choose m+1 } + {d-1 \choose m+1} - {d-s \choose m+1 }  \bigg]\nonumber\\
			 & \quad\quad\quad\quad \quad \quad\quad + \bigg[ \sum_{i=1}^{d-s+2} {d+1-i \choose m} + {d-2 \choose m-1} \bigg] \quad (\text{by \cite[Thm. 3.2]{myLBTpaper}.})\nonumber\\
			&=  {d+1 \choose m+1 } + {d \choose m+1}  +  {d-1 \choose m} + \sum_{i=4}^{d-s+2} {d+1-i \choose m} - {d-s \choose m+1 } \nonumber\\
			&= f_m(\nabla) + \sum_{i=4}^{d-s+2} {d+1-i \choose m} - \sum_{i=m}^{d-s-1} {i \choose m} \nonumber\\
			&= f_m(\nabla) + \sum_{i=s-1}^{d-3} {i \choose m} - \sum_{i=m}^{d-s-1} {i \choose m} \\
			&\geq f_m(\nabla)+ \sum_{i=(s-1) - (s-2)}^{(d-3) - (s-2)} {i \choose m} - \sum_{i=m}^{d-s-1} {i \choose m}\nonumber \\
			&= f_m(\nabla) + \underbrace{\sum_{i=1}^{d-s-1} {i \choose m} - \sum_{i=m}^{d-s-1} {i \choose m} }_{= 0 \; (\text{since } m\geq 1)}\nonumber \\
			&= f_m(\nabla). \nonumber
			\end{align}
			
			\item[Case 5: ] Every facet that does not contain a nonsimple vertex has $d$ vertices.
			
			Let $F_1$ be a facet that does not contain a nonsimple vertex $v$, and let $\mathcal{S}_1 = \{x_0 =v,\; x_1,x_2,\dots  \}$ be the set of vertices outside of $F_1$. By our assumption, $|\mathcal{S}_1|\geq d+1$. If there exists some other facet $F_2$ such that $0<|F_2 \cap \mathcal{S}_1|\leq |\mathcal{S}_1|-2$ (in other words, $F_2$ contains some $x_i$'s but also excludes at least two vertices from $\mathcal{S}_1$.) Then once again by Proposition \ref{thm: key prop}(iii) we obtain $f_m(P)\geq f(\nabla)$.
			
			Assuming such a facet $F_2$ does not exist, the following must be true: \begin{itemize}
				\item[(1)] $|\mathcal{S}_1| = d+1$; 
				\item[(2)] $f_0(P)=2d+1$;
				\item[(3)] there is a unique facet $F_2$ that contains some $x_i$ but not $v$ and $F_2 = \{x_1,\dots ,x_d \}$. 
			\end{itemize} Let $\mathcal{S}_2$ be the set of vertices not in $F_2$ (notice that $\mathcal{S}_2 = V(F_1)\cup v$  and so $|\mathcal{S}_2|=d+1$). Besides $F_1$ and $F_2$, every other facet of $P$ must contain the vertex $v$, and so every other facet contains some vertices from $\mathcal{S}_1$ and $\mathcal{S}_2$. Let $F_3$ be one of them. Then $F_3$ needs to contains at least $d$ vertices (hence exactly $d$ vertices) from each of  $\mathcal{S}_1$ and $\mathcal{S}_2$, so $|F_3| = 2d-1$ and there are exactly two vertices of $P$ outside of $F_3$. This subcase is checked previously in Case 3.
			
		\end{itemize}
		This completes the proof.
	\end{proof}

In summary, we now have the full picture of the lower bounds on face numbers of normal CW $(d-1)$-pseudomanifolds with at least $2d+1$ vertices. The picture is split into two parts by restrictions on the number of facets: \begin{itemize}
	\item[(i)] If the number of facets is $d+2$, by Lemma \ref{lem: bound for (2d+1, d+2)}, the face vector of the polytope, $(T^{d,\frac{d}{2}-2}_2)^*$, is componentwise minimal (when it exists). We know this by completely understanding all possible polytopes (hence all pseudomanifolds by Theorem \ref{thm: spheres are polytopes}) with $d+2$ facets and the fact that all polytopes in this set are completely ordered by $f$-numbers (by Theorem \ref{thm: remark}).
	\item[(ii)] If the number of facets is greater than $d+2$, then by Theorem \ref{thm: 2d+1}, the face vector of a polytope with $d+3$ facets,  $\nabla=(\Stack(T^{d,d-2}_1))^*$, is componentwise minimal.
\end{itemize} 

While the results (i) and (ii) are stated and proved for normal pseudomanifolds with {\bf at least} $2d+1$ vertices, they are sharp only for normal pseudomanifolds with {\bf exactly} $2d+1$ vertices. The reason we need the ``at least'' $2d+1$ part in the statements is to enable us to apply induction. The problem of obtaining sharp lower bounds on the face vectors of $d$-polytopes with $2d+2$ or more vertices is wide open at present: there is not even a plausible conjecture. For polytopes with exactly $2d+2$ vertices, the only result known is the minimal number of edges proved by Pineda-Villavicencio, Ugon, and Yost in 2020 (see \cite{https://doi.org/10.48550/arxiv.2005.06746}).

\section{Concluding remarks}
The most natural open problem is to find the lower bounds for the face numbers of polytopes, strongly regular CW spheres, or even strongly regular normal pseudomanifolds with $2d+2$ and more vertices. In view of our results for $2d+1$ vertices, one can conjecture that the bounds in the case of $2d+2$ vertices are also given by the face numbers of certain polytopes. An approach that is similar to the proof of Theorem \ref{thm: 2d+1} might be feasible but would probably be lengthy and contain quite a few cases.

Another question is the following. By proving Gr\"unbaum's conjecture in \cite{myLBTpaper} we showed the existence of the componentwise minimal $f$-vectors of $d$-polytopes with at most $2d$ vertices. We now know that such an $f$-vector does not exist (in general) when there are $2d+1$ vertices, since the lower bound has two parts. Notice that if in addition we also restrict the number of facets, then the componentwise minimal $f$-vectors still exist. Whether restrictions on the number of facets helps reducing the ``complexity" of the lower bounds for $f$-vectors would be an interesting problem.  This might be more plausible when the number of vertices is relatively small (up to $3d-2$).

\section*{Acknowledgments}
The author would like to thank Isabella Novik for having taken many hours to discuss this project and to help revise the first few drafts of this paper during covid time. The possible generalization of the results from polytopes to spheres was suggested by Jesus De Loera. The author is also grateful to Bennet Goeckner and Marge Bayer for numerous comments and suggestions on the draft, and to Guillermo Pineda-Villavicencio for taking the time to look at the previous version of the paper and to provide feedback.

	\bibliography{lbt_2021}

	\bibliographystyle{acm}

	\newpage
	
	\appendix
	
	\section{\\ }
In this section we will show that the following purely combinatorial inequality (see Proposition \ref{prop: append}) holds. Let $d\geq 3$ and $1\leq m\leq \lceil \frac{d}{2} \rceil$-2, define the following two functions $A(m,d)$ and $B(m,d)$.
\begin{align}
A(m,d) = {d+1 \choose m+1} + {d \choose m+1} + {d-1 \choose m+1} - {\ceil{\frac{d}{2}} \choose m+1} - {\ceil{\frac{d}{2}} -1 \choose m+1}, \text{ as in Lemma \ref{lem: bound for (2d+1, d+2)}.} 
\end{align}

\begin{align}
B(m,d) = {d+1 \choose m+1} + {d \choose m+1} + {d-1 \choose m}, \text{ as in Theorem \ref{thm: 2d+1}.}
\end{align}

We will show that $B(m,d)$ is smaller than $A(m,d)$ roughly for the first half of values of $m$.
\begin{proposition}\label{prop: append}
When $1\leq m < \lceil \frac{d}{2} \rceil-1$ and $d\geq 6$, $B(m,d) \leq A(m,d)$.
\end{proposition}

\begin{proof}
	We will prove the statement by induction on $d$. First let
	\begin{align}
	\delta(m,d) =A(m,d) - B(m,d) =  {d-1 \choose m+1} - {d-1 \choose m} - {\ceil{\frac{d}{2}} \choose m+1} - {\ceil{\frac{d}{2}} -1 \choose m+1}.
	\end{align}
	
When $d=6$ and $m=1$, $\delta(m,d)=1\geq 0$.

Now let $d\geq 7$, and asume that $\delta(m,d-1)\geq 0$ for all $m< \lceil \frac{d-1}{2} \rceil-1$. Consider $m< \lceil \frac{d}{2} \rceil-1$. By Pascal's identity, $\delta(m,d) = \delta(m,d-1) + \delta(m-1,d-1)$. If both $m-1, m < \lceil \frac{d-1}{2} \rceil-1$, then by the inductive hypothesis $\delta(m,d) \geq 0$. The only case that is not covered in the inductive hypothesis is the following.
\begin{align}\left\lceil \frac{d-1}{2} \right\rceil-1=m < \left\lceil \frac{d}{2} \right\rceil-1.
\end{align}
It occurs ony when $d$ is odd (write it as $d= 2a+1$) and $m = \lceil \frac{2a}{2} \rceil-1 = a-1$. Note that $a\geq 4$ since $d\geq 7$. In this case,
\begin{align} \label{eqn: delta}
\delta(a-1,2a+1) &= {2a \choose a} - {2a \choose a-1}  - {a+1 \choose a} - {a \choose a} \nonumber \\
&\stackrel{(*)}{=} \frac{1}{2a+1}{2a+1 \choose a}  - (a+2) \nonumber \\
&= \bigg[ \frac{(2a)(2a-1) \dots (a+3)}{a!} -1 \bigg]\cdot (a+2).
\end{align}

where $(*)$ is by the binomial identity ${\binom {n}{k}}-{\binom {n}{k-1}}={\frac {n+1-2k}{n+1}}{\binom {n+1}{k}}.$

Since $a\geq 4$, 
\begin{align}
\frac{(2a)(2a-1) \dots (a+3)}{a!}  &= 
\frac{2a}{a} \cdot \frac{2a-1}{a-1} \cdot \dots \cdot \frac{a+3}{3}  \cdot \frac{1}{2}\\
&= \frac{2a-1}{a-1} \dots \frac{a+3}{3}  \\
&\geq 1 \quad \text{since every term is at least }1.
\end{align}

Plugging this into (\ref{eqn: delta}), we obtain that $\delta(a-1,2a+1) \geq 0$, which completes the proof.

\end{proof}

\end{document}